\documentclass[12pt]{gen-p-l}
\usepackage{amssymb}
\title[Non-rigid families of Calabi-Yau 3-folds]
{Geometry and arithmetic of non-rigid families of Calabi-Yau 3-folds;
questions and examples}
\author[Eckart Viehweg]{Eckart Viehweg}
\address{Universit\"at Duisburg-Essen, Mathematik, 45117 Essen, Germany}
\email{viehweg@uni-essen.de}
\thanks{This work has been supported by the ``DFG-Schwerpunktprogramm
Globale Methoden in der Komplexen Geometrie'', and by the DFG-Leibniz program.\\[.2cm]
MSC-class: 14D07 (Primary), 14D05, 14J32, 32Q25 (Secondary)}
\author[Kang Zuo]{Kang Zuo}
\address{Universit\"at Mainz,
Fachbereich 17, Mathematik,
55099 Mainz, Germany}
\email{kzuo@mathematik.uni-mainz.de}
\begin{document}
\theoremstyle{plain}
\newtheorem{thm}{Theorem}[]
\newtheorem{theorem}[thm]{Theorem}
\newtheorem{lemma}[thm]{Lemma}
\newtheorem{corollary}[thm]{Corollary}
\newtheorem{proposition}[thm]{Proposition}
\newtheorem{addendum}[thm]{Addendum}
\newtheorem{variant}[thm]{Variant}
\theoremstyle{definition}
\newtheorem{construction}[thm]{Construction}
\newtheorem{notations}[thm]{Notations}
\newtheorem{question}[thm]{Question}
\newtheorem{problem}[thm]{Problem}
\newtheorem{remark}[thm]{Remark}
\newtheorem{remarks}[thm]{Remarks}
\newtheorem{definition}[thm]{Definition}
\newtheorem{claim}[thm]{Claim}
\newtheorem{assumption}[thm]{Assumption}
\newtheorem{assumptions}[thm]{Assumptions}
\newtheorem{properties}[thm]{Properties}
\newtheorem{example}[thm]{Example}
\newtheorem{fact}[thm]{Fact}
\catcode`\@=11
\def\opn#1#2{\def#1{\mathop{\kern0pt\fam0#2}\nolimits}}
\def\bold#1{{\bf #1}}%
\def\underrightarrow{\mathpalette\underrightarrow@}
\def\underrightarrow@#1#2{\vtop{\ialign{$##$\cr
 \hfil#1#2\hfil\cr\noalign{\nointerlineskip}%
 #1{-}\mkern-6mu\cleaders\hbox{$#1\mkern-2mu{-}\mkern-2mu$}\hfill
 \mkern-6mu{\to}\cr}}}
\let\underarrow\underrightarrow
\def\underleftarrow{\mathpalette\underleftarrow@}
\def\underleftarrow@#1#2{\vtop{\ialign{$##$\cr
 \hfil#1#2\hfil\cr\noalign{\nointerlineskip}#1{\leftarrow}\mkern-6mu
 \cleaders\hbox{$#1\mkern-2mu{-}\mkern-2mu$}\hfill
 \mkern-6mu{-}\cr}}}
\let\amp@rs@nd@\relax
\newdimen\ex@
\ex@.2326ex
\newdimen\bigaw@
\newdimen\minaw@
\minaw@16.08739\ex@
\newdimen\minCDaw@
\minCDaw@2.5pc
\newif\ifCD@
\def\minCDarrowwidth#1{\minCDaw@#1}
\newenvironment{CD}{\@CD}{\@endCD}
\def\@CD{\def\A##1A##2A{\llap{$\vcenter{\hbox
 {$\scriptstyle##1$}}$}\Big\uparrow\rlap{$\vcenter{\hbox{%
$\scriptstyle##2$}}$}&&}%
\def\V##1V##2V{\llap{$\vcenter{\hbox
 {$\scriptstyle##1$}}$}\Big\downarrow\rlap{$\vcenter{\hbox{%
$\scriptstyle##2$}}$}&&}%
\def\={&\hskip.5em\mathrel
 {\vbox{\hrule width\minCDaw@\vskip3\ex@\hrule width
 \minCDaw@}}\hskip.5em&}%
\def\verteq{\Big\Vert&&}%
\def\noarr{&&}%
\def\vspace##1{\noalign{\vskip##1\relax}}\relax\let\amp@rs@nd@&\iffalse}\fi
 \CD@true\vcenter\bgroup\relax\let\\=\cr\iffalse}\fi\tabskip\z@skip\baselineskip20\ex@
 \lineskip3\ex@\lineskiplimit3\ex@\halign\bgroup
 &\hfill$\m@th##$\hfill\cr}
\def\@endCD{\cr\egroup\egroup}
\def\>#1>#2>{\amp@rs@nd@\setbox\z@\hbox{$\scriptstyle
 \;{#1}\;\;$}\setbox\@ne\hbox{$\scriptstyle\;{#2}\;\;$}\setbox\tw@
 \hbox{$#2$}\ifCD@
 \global\bigaw@\minCDaw@\else\global\bigaw@\minaw@\fi
 \ifdim\wd\z@>\bigaw@\global\bigaw@\wd\z@\fi
 \ifdim\wd\@ne>\bigaw@\global\bigaw@\wd\@ne\fi
 \ifCD@\hskip.5em\fi
 \ifdim\wd\tw@>\z@
 \mathrel{\mathop{\hbox to\bigaw@{\rightarrowfill}}\limits^{#1}_{#2}}\else
 \mathrel{\mathop{\hbox to\bigaw@{\rightarrowfill}}\limits^{#1}}\fi
 \ifCD@\hskip.5em\fi\amp@rs@nd@}
\def\<#1<#2<{\amp@rs@nd@\setbox\z@\hbox{$\scriptstyle
 \;\;{#1}\;$}\setbox\@ne\hbox{$\scriptstyle\;\;{#2}\;$}\setbox\tw@
 \hbox{$#2$}\ifCD@
 \global\bigaw@\minCDaw@\else\global\bigaw@\minaw@\fi
 \ifdim\wd\z@>\bigaw@\global\bigaw@\wd\z@\fi
 \ifdim\wd\@ne>\bigaw@\global\bigaw@\wd\@ne\fi
 \ifCD@\hskip.5em\fi
 \ifdim\wd\tw@>\z@
 \mathrel{\mathop{\hbox to\bigaw@{\leftarrowfill}}\limits^{#1}_{#2}}\else
 \mathrel{\mathop{\hbox to\bigaw@{\leftarrowfill}}\limits^{#1}}\fi
 \ifCD@\hskip.5em\fi\amp@rs@nd@}
\newenvironment{CDS}{\@CDS}{\@endCDS}
\def\@CDS{\def\A##1A##2A{\llap{$\vcenter{\hbox
 {$\scriptstyle##1$}}$}\Big\uparrow\rlap{$\vcenter{\hbox{%
$\scriptstyle##2$}}$}&}%
\def\V##1V##2V{\llap{$\vcenter{\hbox
 {$\scriptstyle##1$}}$}\Big\downarrow\rlap{$\vcenter{\hbox{%
$\scriptstyle##2$}}$}&}%
\def\={&\hskip.5em\mathrel
 {\vbox{\hrule width\minCDaw@\vskip3\ex@\hrule width
 \minCDaw@}}\hskip.5em&}
\def\verteq{\Big\Vert&}
\def\novarr{&}
\def\noharr{&&}
\def\SE##1E##2E{\slantedarrow(0,18)(4,-3){##1}{##2}&}
\def\SW##1W##2W{\slantedarrow(24,18)(-4,-3){##1}{##2}&}
\def\NE##1E##2E{\slantedarrow(0,0)(4,3){##1}{##2}&}
\def\NW##1W##2W{\slantedarrow(24,0)(-4,3){##1}{##2}&}
\def\slantedarrow(##1)(##2)##3##4{%
\thinlines\unitlength1pt\lower 6.5pt\hbox{\begin{picture}(24,18)%
\put(##1){\vector(##2){24}}%
\put(0,8){$\scriptstyle##3$}%
\put(20,8){$\scriptstyle##4$}%
\end{picture}}}
\def\vspace##1{\noalign{\vskip##1\relax}}\relax\let\amp@rs@nd@&\iffalse}\fi
 \CD@true\vcenter\bgroup\relax\let\\=\cr\iffalse}\fi\tabskip\z@skip\baselineskip20\ex@
 \lineskip3\ex@\lineskiplimit3\ex@\halign\bgroup
 &\hfill$\m@th##$\hfill\cr}
\def\@endCDS{\cr\egroup\egroup}
\newdimen\TriCDarrw@
\newif\ifTriV@
\newenvironment{TriCDV}{\@TriCDV}{\@endTriCD}
\newenvironment{TriCDA}{\@TriCDA}{\@endTriCD}
\def\@TriCDV{\TriV@true\def\TriCDpos@{6}\@TriCD}
\def\@TriCDA{\TriV@false\def\TriCDpos@{10}\@TriCD}
\def\@TriCD#1#2#3#4#5#6{%
\setbox0\hbox{$\ifTriV@#6\else#1\fi$}
\TriCDarrw@=\wd0 \advance\TriCDarrw@ 24pt
\advance\TriCDarrw@ -1em
\def\SE##1E##2E{\slantedarrow(0,18)(2,-3){##1}{##2}&}
\def\SW##1W##2W{\slantedarrow(12,18)(-2,-3){##1}{##2}&}
\def\NE##1E##2E{\slantedarrow(0,0)(2,3){##1}{##2}&}
\def\NW##1W##2W{\slantedarrow(12,0)(-2,3){##1}{##2}&}
\def\slantedarrow(##1)(##2)##3##4{\thinlines\unitlength1pt
\lower 6.5pt\hbox{\begin{picture}(12,18)%
\put(##1){\vector(##2){12}}%
\put(-4,\TriCDpos@){$\scriptstyle##3$}%
\put(12,\TriCDpos@){$\scriptstyle##4$}%
\end{picture}}}
\def\={\mathrel {\vbox{\hrule
   width\TriCDarrw@\vskip3\ex@\hrule width
   \TriCDarrw@}}}
\def\>##1>>{\setbox\z@\hbox{$\scriptstyle
 \;{##1}\;\;$}\global\bigaw@\TriCDarrw@
 \ifdim\wd\z@>\bigaw@\global\bigaw@\wd\z@\fi
 \hskip.5em
 \mathrel{\mathop{\hbox to \TriCDarrw@
{\rightarrowfill}}\limits^{##1}}
 \hskip.5em}
\def\<##1<<{\setbox\z@\hbox{$\scriptstyle
 \;{##1}\;\;$}\global\bigaw@\TriCDarrw@
 \ifdim\wd\z@>\bigaw@\global\bigaw@\wd\z@\fi
 \mathrel{\mathop{\hbox to\bigaw@{\leftarrowfill}}\limits^{##1}}
 }
 \CD@true\vcenter\bgroup\relax\let\\=\cr\iffalse}\fi
 \tabskip\z@skip\baselineskip20\ex@
 \lineskip3\ex@\lineskiplimit3\ex@
 \ifTriV@
 \halign\bgroup
 &\hfill$\m@th##$\hfill\cr
#1&\multispan3\hfill$#2$\hfill&#3\\
&#4&#5\\
&&#6\cr\egroup%
\else
 \halign\bgroup
 &\hfill$\m@th##$\hfill\cr
&&#1\\%
&#2&#3\\
#4&\multispan3\hfill$#5$\hfill&#6\cr\egroup
\fi}
\def\@endTriCD{\egroup}
\newcommand{\sA}{{\mathcal A}}
\newcommand{\sB}{{\mathcal B}}
\newcommand{\sC}{{\mathcal C}}
\newcommand{\sD}{{\mathcal D}}
\newcommand{\sE}{{\mathcal E}}
\newcommand{\sF}{{\mathcal F}}
\newcommand{\sG}{{\mathcal G}}
\newcommand{\sH}{{\mathcal H}}
\newcommand{\sI}{{\mathcal I}}
\newcommand{\sJ}{{\mathcal J}}
\newcommand{\sK}{{\mathcal K}}
\newcommand{\sL}{{\mathcal L}}
\newcommand{\sM}{{\mathcal M}}
\newcommand{\sN}{{\mathcal N}}
\newcommand{\sO}{{\mathcal O}}
\newcommand{\sP}{{\mathcal P}}
\newcommand{\sQ}{{\mathcal Q}}
\newcommand{\sR}{{\mathcal R}}
\newcommand{\sS}{{\mathcal S}}
\newcommand{\sT}{{\mathcal T}}
\newcommand{\sU}{{\mathcal U}}
\newcommand{\sV}{{\mathcal V}}
\newcommand{\sW}{{\mathcal W}}
\newcommand{\sX}{{\mathcal X}}
\newcommand{\sY}{{\mathcal Y}}
\newcommand{\sZ}{{\mathcal Z}}
\newcommand{\A}{{\mathbb A}}
\newcommand{\B}{{\mathbb B}}
\newcommand{\C}{{\mathbb C}}
\newcommand{\D}{{\mathbb D}}
\newcommand{\E}{{\mathbb E}}
\newcommand{\F}{{\mathbb F}}
\newcommand{\G}{{\mathbb G}}
\newcommand{\HH}{{\mathbb H}}
\newcommand{\I}{{\mathbb I}}
\newcommand{\J}{{\mathbb J}}
\renewcommand{\L}{{\mathbb L}}
\newcommand{\M}{{\mathbb M}}
\newcommand{\N}{{\mathbb N}}
\newcommand{\BP}{{\mathbb P}}
\newcommand{\Q}{{\mathbb Q}}
\newcommand{\R}{{\mathbb R}}
\newcommand{\T}{{\mathbb T}}
\newcommand{\U}{{\mathbb U}}
\newcommand{\V}{{\mathbb V}}
\newcommand{\W}{{\mathbb W}}
\newcommand{\X}{{\mathbb X}}
\newcommand{\Y}{{\mathbb Y}}
\newcommand{\Z}{{\mathbb Z}}
\newcommand{\id}{{\rm id}}
\newcommand{\rk}{{\rm rank}}
\newcommand{\END}{{\mathbb E}{\rm nd}}
\newcommand{\End}{{\rm End}}
\newcommand{\Hg}{{\rm Hg}}
\newcommand{\tr}{{\rm tr}}
\newcommand{\Sl}{{\rm Sl}}
\newcommand{\Gl}{{\rm Gl}}
\newcommand{\Cor}{{\rm Cor}}
\newcommand{\Hom}{{\sH}{\rm om}}
\newcommand{\s}{{\rm sl}}
\newcommand{\ch}{{\rm c}}
\begin{abstract}
We speculate about the structure of maximal
product subvarieties of moduli stacks of Calabi-Yau manifolds.
We discuss an example of a family of quintic hypersurfaces in $\BP^4$,
parameterized by the product of two ball quotients, one of dimension two,
the second one of dimension $12$.
\end{abstract}
\maketitle
Let $\sM_h(\C)$ denote the set of isomorphism classes of minimal
polarized manifolds $F$ with fixed Hilbert polynomial $h$,
and let $\sM_h$ be the corresponding moduli functor, i.e.
$$
\sM_h(U)=\left\{ \begin{array}{c}(f:V\to U,\sL) ; f \mbox{ smooth and}\\
(f^{-1}(u),\sL|_{f^{-1}}(u))\in \sM_h(\C), \mbox{ for all }
u \in U\end{array}\right\}
$$
There exists a quasi-projective coarse moduli scheme
$M_h$ for $\sM_h$. Fixing a projective manifold $\bar{U}$ and
the complement $U$ of a normal crossing divisor, we want to consider
$${\rm\bf H}=\left\{ \begin{array}{c}
\varphi:(\bar U, U) \to (\overline M_h, M_h) \quad{\rm induced }
\\
{\rm by\, polarized\, families}\, f:X \to U
\end{array} \right\}.$$
Since $M_h$ is just a coarse moduli scheme, it is not clear whether
${\rm\bf H}$ has a scheme structure. However, by \cite{Popp}, if all
$F\in \sM(\C)$ admit a locally injective Torelli map, there exists
a fine moduli scheme $M_h^{N}$ with a level structure $N$ and \'etale
over $M_h.$ By abuse of notations, we will replace $\sM_h$ by the moduli
functor of polarized manifolds with a level $N$ structure, and fix some
compactification $\overline M_h$. Then
${\rm\bf H}$ parameterizes all morphisms from
$$
\varphi:(\bar U, U) \>>> (\overline M_h, M_h),
$$
hence it is a scheme. Moreover there exists a universal family $f:X\to {\rm\bf H}\times U.$

As Kov\'acs, Bedulev-Viehweg, Oguiso-Viehweg, and Viehweg-Zuo have shown
${\rm\bf H}$ is of finite type.

\begin{definition}\label{def}
 $\varphi: U\to M_h$ called rigid if the
component of ${\rm\bf H}$ containing $\varphi$ is
zero-dimensional.
\end{definition}
\begin{problem}\label{quest}
Study the geometry of ${\rm\bf H}$ and the arithmetic properties
(for example the Mumford-Tate group)  of the universal family
$f:X\>>> {\rm\bf H}\times U.$
\end{problem}
{\bf Acknowledgment.} We thank the referee for pointing out several misprints and ambiguities.

This note was finished when the second
named author visited the Institute of Mathematical Sciences at the
Chinese University of Hong Kong. He would like to thank the
members of the Institute for their hospitality.

\section{Splitting of variations of Hodge structures}

Let us start by recalling some of the properties of
complex polarized variations of Hodge structures, and of families
of Calabi-Yau manifolds.

\begin{proposition}\label{prop1} If $\V$ is an irreducible complex
polarized variation of Hodge structures over $U_1\times
\cdots\times U_\ell$ then
$$\V=p^*_1(\V_1)\otimes\cdots\otimes p^*_\ell(\V_\ell),$$
for complex polarized variations of Hodge structures
$\V_i$ over $U_i.$
\end{proposition}

\begin{proof} The proof (see \cite{VZ}, Section 3, for the details) uses Schur's Lemma and
Deligne's semi-simplicity of complex polarized variations of Hodge structures.
\end{proof}
\section{Products in moduli stacks of Calabi-Yau manifolds}

Since Calabi-Yau manifolds are un-obstructed, the fine moduli scheme
$M_h$ is smooth, and we choose a smooth projective compactification
$\overline M_h$ such that $\overline M_h \setminus M_h$ is a normal crossing divisor.
Let $g:\sX\to M_h$ be the universal family. We will assume moreover, that the local
monodromies of $R^mg_*\C_{\sX}$ around the components of $\overline M_h \setminus M_h$
are uni-potent, where $m=\deg(h)$ is the dimension of the fibres.

Consider a smooth family 
$$
f:X\>>> U_1\times \cdots\times U_\ell=U
$$ 
of Calabi-Yau $m$-folds, such that $\varphi: U\to M_h$ is
generically finite. We assume that the factors $U_i$ are non singular, and that $\dim(U_i)>0$. Let
$\V\subset R^mf_*(\C_X)$ be the irreducible  sub variation of Hodge
structures with system of Hodge bundles $$\bigoplus_{p+q=m}E^{p,q}$$
such that  $E^{m,0}=f_*\Omega^m_{X/U}.$
\begin{fact}
The Kodaira-Spencer map is injective and factors through
$$d\varphi: T_U\>>>  E^{m-1,1}\otimes {E^{m,0}}^{-1}\subset
\varphi^*T_{M_h}.$$
\end{fact}
By Proposition \ref{prop1} one has a decomposition
$\V=\V_1\otimes\cdots\otimes\V_\ell$. Let us write
$$\bigoplus_{p+q=m}F_i^{p,q}$$ for the system of Hodge bundles
of $\V_i$, and $\varphi_i:U \to U_i \to \sD_i$ for the corresponding
period map. Then
$$d\varphi_i: T_{U_i}\>>>  F^{m_i-1,1}_i\otimes {F^{m_i,0}_i}^{-1}\subset \varphi^*_iT_{\sD_i},
$$
for $1\leq i\leq \ell$. As in \cite{VZ}, 3.5, a comparison of Hodge bundles on both sides gives rise
to
\begin{proposition}\label{prop2} \ 
\begin{enumerate}
\item[ i.] The cup-product
$$
\bigoplus_{1\leq i_1<\cdots< i_k\leq \ell}
T_{U_{i_1}}\otimes\cdots\otimes T_{U_{i_k}} \>>> R^k f_*T^k_{X/U}
$$
is injective for $1\leq k\leq \ell.$
\item[ii.] If $\varphi:U_1\times \cdots\times U_\ell\to M_h$ is
an embedding and if $\ell=m$ is the dimension of the fibres of $f$ then $U_1\times \cdots\times
U_\ell$ is a product of curves, and uniformized by
$\V_1\oplus\cdots\oplus\V_\ell$ over an algebraic number field.
\end{enumerate}
\end{proposition}

\begin{question}
When will $U_1\times \cdots\times U_m$ be a
product of Shimura curves?
\end{question}
\begin{remark}
A similar argument shows that part i) of Proposition \ref{prop2} also holds true
for moduli stacks of hyper-surfaces in $\BP^n$ (see \cite{VZ}, 3.5 b)).
\end{remark}
\begin{question}
Does Proposition \ref{prop2}, 1) hold true for moduli
stacks of minimal polarized manifolds?
\end{question}
If $U_1\times \cdots\times U_\ell$ maps generically finite
to a moduli stack $\sM_h$ of minimal polarized manifolds,
then it has been shown in \cite{VZ0}, Corollary 6.4, that 
$$\ell \leq m=\deg(h).$$
\begin{question}
Can one improve this bound for certain moduli stacks
and, fixing $\ell$, what are optimal bounds for the dimensions of the
$U_i$?
\end{question}
Since we assumed $M_h$ to be a fine moduli space, 
deformations of the morphism $\varphi:U \to M_h$ correspond to deformations of the
family $f: X \to U$. If one assumes that $U$ has a compactification
$\bar U$ such that $\varphi$ extends to $\varphi: \bar U \to \overline M_h$,
in such a way that the pre-image of $S=\overline M_h \setminus M_h$
remains a reduced normal crossing divisor, the first order deformations of the first type
are classified by $$H^0(\bar U,\varphi^*T_{\overline M_h}(-\log S)).$$

\begin{proposition}\label{prop3}
Assume in addition that $f$ extends to a proper morphism
$f:\bar X\to\bar U$, semi-stable in codimension one, and that $f^*f_*\omega_{\bar X/\bar
U}\to \omega_{\bar X/\bar U}$ is an isomorphism outside of
$f^{-1}(Z)$ for some $Z\subset \bar U$ closed and of codimension
at least two. Then
$$ \dim H^0(\bar U,\varphi^*T_{\overline M_h}(-\log S))$$
is invariant  under infinitesimal deformations.

In  particular, by Ran's $T^1$-lifting property,
deformations of those families $f: X \to U$ of  Calabi-Yau manifolds
with $U$ fixed are un-obstructed.
\end{proposition}
\begin{remark}
We expect that Proposition \ref{prop3} holds true under
weaker and more natural conditions on the boundary.
\end{remark}
\begin{proof} Since we are only interested in global sections,
taking complete intersection we may assume that $\dim\bar
U=1$, that all fibres are semi-stable and that
$$f^*f_*\omega_{\bar X/\bar U}\>>> \omega_{\bar X/\bar U}$$ is an
isomorphism.

Recall that (choosing a level $N$ structure) we assumed the existence of a
universal family $f: \sX \to M_h$. The pull back of the logarithmic Higgs
field
$$\theta: E\>>> E\otimes\Omega^1_{\overline M_h}(\log S)$$ of the variation of Hodge structures
$R^mf_*\Q_\sX$ to $\bar U$ corresponds to  a sub-sheaf
$$ \varphi^*T_{\overline M_h}(-\log S)\>>> (\mathcal End(\varphi^* E),\theta^{\mathcal End}).$$
By (\cite{Z}, Prop. 2.1) $\theta^{\mathcal End}(\varphi^*T_{\overline
M_h}(-\log S))=0.$ This means that the above sub-sheaf is a Higgs sub-sheaf.

We need the following theorem on intersection cohomology and Higgs
cohomology of a complex polarized variation of Hodge structures $\W$ with uni-potent local
monodromy around $S$. Let $(F,\theta)$ denote the logarithmic
Higgs bundle of $\W.$ We consider the complex of sheaves defined
by the Higgs field
$$ F\>{\theta}>> F\otimes\Omega^1_{\bar U}(\log S)\>{\theta}>>
F\otimes \Omega^2_{\bar U}(\log S)\>>>\cdots .$$ In \cite{Zu} (for
$\dim\bar U=1$ in an implicit way) and in \cite{JYZ} (in general)
one finds the definition of an algebraic $L_2$- sub complex of
sheaves

$$\begin{CDS} F\>{\theta}>> F\otimes\Omega^1_{\bar U}(\log S)\>{\theta}>> F\otimes \Omega^2_{\bar U}(\log
S)\>>>\cdots\\
\cup\novarr \novarr  \cup\novarr \novarr\cup\\
F_{(2)}\>{\theta}>> (F\otimes\Omega^1_{\bar U}(\log
S))_{(2)}\>{\theta}>> (F\otimes \Omega^2_{\bar U}(\log
S))_{(2)}\>>>\cdots
\end{CDS}$$
determined by an algebraic  condition on $F|_S$ imposed by the weight-filtration of
$${\rm res}(\theta): F|_S\>>> \varphi^* E|_S.$$  Note that for a
sub sheaf  $F'\subset {\rm Ker}(\theta)$, one has $F'\subset
F_{(2)}.$
\begin{theorem}[\cite{Zu} for $\dim \bar U=1$, \cite{JYZ}]\label{theorem4}
$$  \mathbb H^i\big(F_{(2)}\>{\theta}>> (F\otimes\Omega^1_{\bar U}(\log
S))_{(2)}\>{\theta}>>\cdots\big)\simeq H^i_{\rm intersection}(\W).$$
\end{theorem}

 Back to our situation,  the exact sequence of complexes of sheaves
 $$0\>>> (\varphi^*T_{\overline M_h}(-\log S),0)\>>> (\mathcal End(\varphi^* E),\theta^{\mathcal End})\>>> (Q,\theta)\>>> 0$$
gives rise to a long  exact sequence

$$\begin{array}{cccccc}
 \cdots &\to & \mathbb H^{i-1}(Q,\theta)&\to& H^i(\varphi^*T_{\overline M_h}(-\log S))&\to \mathbb
 H^i(\mathcal End (\varphi^* E)_{(2)},\theta^{\mathcal End})\\
         &\to &\mathbb H^i(Q,\theta)    &\to &   H^{i+1}(\varphi^*T_{\overline M_h}(-\log S))&\to \mathbb
 H^{i+1}(\mathcal End (\varphi^* E)_{(2)},\theta^{\mathcal End})\\
&\to &\cdots &&&
\end{array}$$
Since we assumed the fibres $f^{-1}(p)$ of $f$ to be semi-stable and minimal,
\cite{KN} implies that $f^{-1}(p)$ has no obstruction to deformations in any
direction. This means that the pullback of the Kodaira-Spencer
map of the moduli space to $\bar U$
$$(\varphi^*T_{\overline M_h}(-\log S),0)\>>> \mathcal
(\mathcal End(\varphi^* E)_{(2)},\theta^{\mathcal End})\>>> (\varphi^*E^{m-1,1}\otimes
\varphi^*E^{0,m},0)$$ is an isomorphism. Taking in account that those are maps between
complexes of sheaves, we find
$$ H^i(\varphi^*T_{\overline M_h}(-\log S))\>>> \mathbb H^i(\mathcal End (\varphi^*E)_{(2)},
 \theta^{\mathcal End})$$
to be  injective for all $i$. Hence there is a  splitting

$$\mathbb H^i(\mathcal End (\varphi^*E)_{(2)},
 \theta^{\mathcal End})=H^i((\varphi^*T_{\overline M_h}(-\log S))\oplus
 \mathbb H^i(Q,\theta).$$

 By Theorem \ref{theorem4} $\mathbb H^i(\mathcal End (\varphi^*E)_{(2)},
 \theta^{\mathcal End})$ is isomorphic to the intersection
 cohomology, hence is invariant under small deformations.
 Using the semi continuity of the hyper-cohomology of complexes of
 sheaves one shows that both  $H^i((\varphi^*T_{\overline M_h}(-\log S))$
 and  $\mathbb H^i(Q,\theta)$ are invariant under small
 deformations.
\end{proof}
\begin{corollary}\label{cor4}
Under the assumptions made in \ref{prop3} the scheme $\bf \rm H$ is smooth.
\end{corollary}
\section{Applications}

Again $f: X\to U$ denotes a smooth family of Calabi-Yau $3$-folds, such that $\varphi: U\to
M_h$ is generically finite. We keep the assumption, that $M_h$ has a universal family.
Moreover, we choose a compactification $\overline M_h$ with $\overline M_h\setminus M_h$ a normal
crossing divisor, such that $U\to M_h$ extends to $\bar U \to \overline M_h$.

Staring with
$${\rm\bf H_1}={\rm Hom} ((\overline U,U), (\overline M_h,
M_h)),$$
consider
$${\rm\bf H_2}={\rm Hom} (({\rm\bf \overline H_1}\times\{0\},{\rm\bf H_1}\times\{0\}),
(\overline M_h, M_h)),\quad \{0\}\in U,$$
together with the induced family $f: X\to {\rm \bf H_1}\times{\rm\bf H_2=H}.$

Let $\V\subset R^3f_*(\C_X)$ be the irreducible sub variation of Hodge structures
with Hodge decomposition
$$\bigoplus_{p+q=3}F^{p,q}\quad
{\rm with}\quad  F^{3,0}=f_*\Omega^3_{X/{\rm\bf H}}.$$

Recall that by Proposition \ref{prop1} one has a decomposition $\V=\V_1\otimes\V_2$, where
$\V_i$ is the pull back of a $\C$ variation of Hodge structures on ${\rm\bf H_i}$.
Comparing the possible Hodge numbers, one finds:

\begin{proposition}\label{prop4}
$\V_i$ has one of the following  Hodge types:
\begin{enumerate}
\item[a.] $F_i^{1,0}\oplus F_i^{0,1},\quad {\rm rk} F_i^{1,0}=1.$
\item[b.] $F_i^{2,0}\oplus F_i^{1,1}\oplus
F_i^{0,2},\quad {\rm rk} F_i^{2,0}={\rm rk} F_i^{0,2}=1$, and $\V_i$ is real.
\item[c.] Moreover, if $\V_1$ is of type b), then ${\rm  rk}\V_2=2.$
\end{enumerate}
\end{proposition}
It is well known that the period domains $\sD_i$ of Hodge structures of types a) or b)
are the bounded symmetric domain of the algebraic group $U(1,{\rm rk} \V^{0,1}_i)$,
or $SO(2,{\rm rk} \V^{1,1}_i)$, respectively.

The un-obstructedness for deformations of families implies that the generically finite period map
${\rm\bf \tilde H_i}\to \sD_i$ has to be dominant. Let us assume that $U\to M_h$ is injective.
\begin{question}\label{conj1} \
\begin{enumerate}
\item[1.] Is ${\rm\bf H_i^s}\simeq \sD_i/\Gamma_i$ for some $\Gamma_i$ a partial compactification of
${\rm\bf H_i}$?
\item[2.] What is the moduli-interpretation of points in ${\rm\bf H^s_i}\setminus {\rm\bf H_i}?$
\end{enumerate}
\end{question}
\section{An example of a non-rigid family of Calabi-Yau quintic threefolds}

Let $f_5(x_2,x_1,x_0)\in \C[x_2,x_1,x_0]$  be the polynomial of a quintic plane curve in $\BP^2.$
Then $$x_3^5+f_5(x_2,\,x_1,\,x_0)$$ defines a quintic hypersurface in $\BP^3$, and
$$x^5_4+x^5_3+f_5(x_2,\,x_1,\,x_0)$$ a Calabi-Yau quintic $3$-fold in $\BP^4.$

Obviously this construction can also be done locally over the moduli stack $M_{5,2}$
of quintic plane curves in $\BP^2$, starting with the universal family $f: X\to M_{5,2}$
of curves. Replacing $M_{5,2}$ by some covering, one can glue those families
as family of subvarieties in some projective bundle (see \cite{VZ}).
The resulting family of surfaces will be denoted by $g_1:Z_1\to M_{5,2}$, and
the one of threefolds by $g_2: Z_2\to M_{5,2}$.

\begin{remark} As pointed out by S.T. Yau, this family has been studied by S. Ferrara
and J. Louis \cite{FL}. They have shown that the Yukawa-coupling is zero and that
and the monodromy lies in $SU(2,1)$. In \cite{VZ} the exact length of the Yukawa coupling
is calculated for such families.
\end{remark}

One can play a similar game, starting with $5$ points in $\BP^1.$
say with equation $h_5(x_1,\,x_0)\in \C[x_1,\,x_0]$. Then $x_2^5+h_5(x_1,\,x_0)$
defines a quintic plane curve. Again, one can do such a construction starting with the
universal family $P\to M_{5,1}$ of 5 points in $\BP^1$, and one obtains
a family $g_0: Z_0\to M_{5,1}$ of quintic plane curves.

Finally $\Sigma_5$ denotes the Fermat curve $x_2^5+x_1^5+x_0^5=0$ of degree $5$.
\begin{proposition}\label{prop5}
The fibre product $Z_1\times\Sigma_5\to M_{5,2}$ admits an $\Z_5$-action over $M_{5,2}$,
given fibrewise by
$$(x_3, x_2, x_1, x_0), (y_2, y_1, y_0)\mapsto(e^{2\pi i/5}x_3,x_2, x_1, x_0), (e^{2\pi i/5} y_2, y_1, y_0).$$
\begin{enumerate}
\item[1] The family of Calabi-Yau quintics $g_2: Z_2\to
M_{5,2}$ can be reconstructed as:
$$
\begin{CDS}
 (Z_1\times\Sigma_5)/\Z_5 \< {\rm blow up}   << \widehat {(Z_1\times\Sigma_5)/\Z_5}\>{\rm blow down}>> Z_2\\
\novarr \SE E g_1 E  \quad\V V  V       \quad   \SW W g_2 W\\
\novarr \novarr \quad M_{5,2}.\end{CDS}
$$
\item[2.] The  construction in 1) extends to the
product family
$$\begin{CDS}
(Z_1\times Z_0)/\Z_5 \< {\rm blow \,up}   << \widehat {(Z_1\times Z_0)/\Z_5}\>{\rm blow\, down}>> \sZ_2\\
\novarr \SE E (g_1,g_0)  E  \quad\V V  V       \quad   \SW W h_2 W\\
\novarr \novarr \quad
M_{5,2}\times M_{5,1}.\\
\end{CDS}
$$
\item[3.] The family $h_2: \sZ_2\to M_{5,2}\times M_{5,1}$ of Calabi-Yau
quintics is a universal family of the form
$$h_2: \sZ_2\>>> {\rm\bf H_1}\times {\rm\bf H_2},$$
i.e. for suitable compactifications $\overline M_h$, $\overline M_{5,2}$ and $\overline M_{5,1}$
and for some base point $u \in M_{5,2}$ and $u'\in M_{5,1}$
\begin{gather*}
M_{5,2}={\rm\bf H_1}={\rm Hom}((\{u\}\times\overline M_{5,1},\{u\}\times M_{5,1}),
(\overline M_h, M_h)), \\ \mbox{and \ \ }\\
M_{5,1}={\rm\bf H_2}={\rm Hom}((\overline M_{5,2}\times\{u'\},M_{5,2}\times\{u'\}),
(\overline M_h, M_h)).
\end{gather*}
Moreover, a partial compactification ${\rm\bf H_1^s}$ of ${\rm\bf H_1}$ is a
2-dimensional complex arithmetic ball quotient, and a partial compactification ${\rm\bf H_2^s}$
of ${\rm\bf H_2}$ is a 12-dimensional complex non-arithmetic ball quotient.
\end{enumerate}
\end{proposition}
\begin{proof}[Sketch of the proof]
1) and 2) have been shown in (\cite{VZ},
Proposition 6.4). For 3) consider the eigen-space decompositions
$$
R^1g_{0*}(\bar\Q_{Z_0})=\bigoplus_{i=1}^4R^1g_{0*}(\bar\Q_{Z_0})_i,\mbox{ \ \ and \ \ }
R^2g_{1*}(\bar\Q_{Z_1})=\bigoplus_{i=1}^4R^2g_{1*}(\bar\Q_{Z_1})_i
$$
for the $\Z_5$-action. Recall that the restriction of
$R^1g_{0*}(\bar\Q_{Z_0})_i$ to a point in $M_{5,1}$ is a Hodge structure
with $H^0(\BP^1,\Omega^1_{\BP^1}(5-i))$ in degree $(1,0)$ and
$H^1(\BP^1,\sO_{\BP^1}(-i))$ in degree $(0,1)$. Hence
$R^1g_{0*}(\bar\Q_{Z_0})_i$ is unitary for $i=1$ and $i=4$.

$R^1g_{0*}(\bar\Q_{Z_0})_3$ and $R^1g_{0*}(\bar\Q_{Z_0})_2$ are not unitary, and dual to each other.
By Deligne-Mostow \cite{DM} $M_{5,1}^s$ is uniformized by
$R^1g_{0*}(\bar\Q_{Z_0})_3$ as a $2$-dimensional arithmetic ball quotient,
which is a component of the moduli space, parameterizing Abelian
varieties of dimension $6$ with complex multiplication $\Q(\zeta)$ for $\zeta=e^{2\pi i/5}$.

For $g_1$ the situation is more complicated. One easily computes that the Higgs bundle
of $R^2g_{1*}(\bar\Q_{Z_1})_2$ is trivial in degree $(0,2)$, of a rank one in degree $(2,0)$ and of a 
rank $12$ in degree $(1,1)$. 
A similar argument as the one used by Deligne-Mostow allows to show that $M^s_{5,2}$ is
uniformized by $R^2g_{1*}(\bar\Q_{Z_1})_2$ as a $12$-dimensional complex
ball quotient. 

There is a Galois conjugate $R^2g_{1*}(\bar\Q_{Z_1})_2^\sigma$, which is neither the dual of
$R^2g_{1*}(\bar\Q_{Z_1})_2$, nor unitary. As in Deligne-Mostow this implies that
the ball quotient is not arithmetic.

The quotient by $\Z_5$, together with  blowing up and blowing
down, gives rise to a $\Q$-Hodge isometry (see \cite{VZ}, 7.4)
$(R^3h_{2*}\Q_{\sZ_2})\simeq \V'\oplus\T$ with
\begin{gather*}
\V'\otimes \Q(\zeta)= \bigoplus_{i=1}^4 \V_{i,5-i} \mbox{ \ \ for \ \ }
\V_{i,5-i}=R^2g_{1*}(\Q(\zeta)_{Z_1})_i\otimes R^1g_{0*}(\Q(\zeta)_{Z_0})_{5-i} \\
\mbox{ and with \ \ } \T= \bigoplus^4R^1f_*\Q_X(1),
\end{gather*}
where $(1)$ denotes the Tate-twist. Of course, we should write
$$
\V_{i,5-i}={\rm pr}_1^*R^2g_{1*}(\Q(\zeta)_{Z_1})_i\otimes {\rm pr}_2^*R^1g_{0*}(\Q(\zeta)_{Z_0})_{5-i},
$$
but we suppress the pullback under the projections in our notation. Remark that
$\T$ is the part of the variation of Hodge structures, coming from the
blowing ups.
$\V_{i,5-i}$ is an irreducible sub-variation of Hodge structures in $(R^3h_{2*}\Q(\zeta)_{\sZ_2})$,
and for the corresponding $\C$ variation of Hodge structures, the fibre over a point $y$ has
$$
H^0(\BP^2,\Omega^2_{\BP^2}(5-i))\otimes H^0(\BP^1,\Omega^1_{\BP^1}(i))
$$
in degree $(3,0)$. This is zero for $i=1$ and $i=4$, and the $\C$ variation of Hodge structures given by
$$
\V=\V_{2,3}\oplus \bar \V_{2,3}= \V_{2,3}\oplus \V_{3,2}
$$
contains the Hodge bundle ${h_{2*}\Omega^3_{\sZ_2/ M_{5,2}\times M_{5,1}}}$.

By abuse of notations we will regard $\V$ and the $\V_{i,5-i}$ as $\bar\Q$
variations of Hodge structures.
Write $(R^3h_{2*}\bar\Q_{\sZ_2})=\V\oplus\W$, and let
$(F^{2,1}\oplus F^{1,2}\oplus F^{0,3},\theta)$ denote the system
of Hodge bundles corresponding to $\W$. The missing part of \ref{prop5}, 3),
follows from the next two Claims.
\end{proof}
\begin{claim}\label{claim6}
There is no nontrivial extension
$$
\begin{CD}
\sZ_2 \> {h_2}   >>  M_{5,2}\times M_{5,1}\\
\V \subset VV \V V \subset V      \\
\sZ'_2 \> {h'_2}   >>   \quad N\times M_{5,1},
\end{CD}
$$
such that the induced morphism
$\varphi: N\times M_{5,1}\to M_h$ is generically finite over its image.
\end{claim}
\begin{proof}
A deformation $N\times M_{5,1}$ of $M_{5,1}=\{u\}\times M_{5,1}$, which does not lie in
$M_{5,2}\times M_{5,1}$, induces a non-zero flat section $\tau$ of ${\rm
End}(\V\oplus\W)|_{M_{5,1}} $ of type $(-1,1)$, which does not
respect the direct sum decomposition $\V\oplus \W$.

In fact, if it does, one has $\tau(\V|_{M_{5,1}})\subset \V|_{M_{5,1}}$. The restriction of
$\V_{2,3}$ and $\V_{3,2}$ to ${M_{5,1}}$ are direct sums of local systems
isomorphic to $R^1g_{0*}(\bar\Q_{Z_0})_3$ or $R^1g_{0*}(\bar\Q_{Z_0})_2$, respectively.
As uniformizing variations of Hodge structures of a ball quotient, both are irreducible
and, since the Hodge numbers are different, $R^1g_{0*}(\bar\Q_{Z_0})_3$ is not
isomorphic to $R^1g_{0*}(\bar\Q_{Z_0})_2$. So $\tau$ respect the decomposition
$\V=\V_{2,3}\oplus \V_{3,2}$ and it is induced by an endomorphism
of $R^2g_{1*}(\bar\Q_{Z_1})_2|_{\{u\}}$ of type $(-1,1)$. 
The calculation of Hodge numbers, indicated above, shows that those are lying in a $12$ dimensional 
vector space. So they correspond all to the deformations along $M_{5,2}$.

On the other hand, if $\tau$ does not respect the direct sum decomposition, one obtains
a non trivial morphism $\V|_{M_{5,1}} \to \W|_{M_{5,1}}$.
Since $\W|_{M_{5,1}}$ is a direct sum of unitary
local systems, and since for $i=2$ or $i=3$ there exists no non-trivial morphism
$$R^1g_{0*}(\bar\Q_{Z_0})_i\>>> {\rm unitary\, local \, system},$$
this leads to a contradiction.
\end{proof}
\begin{claim}\label{claim7}
There is no nontrivial extension
$$
\begin{CD}
\sZ_2 \> {h_2}   >>  M_{5,2}\times M_{5,1}\\
\V \subset VV \V V \subset V      \\
\sZ'_2 \> {h'_2}   >>   \quad M_{5,2}\times N,
\end{CD}
$$
such that the induced morphism
$\varphi: M_{5,2}\times N\to M_h$ is generically finite over its image.
\end{claim}
\begin{proof}
Again, the deformations of $M_{5,2}\times \{u\}$ correspond to flat section $\tau$ of
${\rm End}(\V\oplus\W)|_{M_{5,2}} $ of type $(-1,1)$. If $\tau$ respects the direct sum decomposition,
the irreducibility of $R^2g_{1*}(\bar\Q_{Z_1})_2$ and $R^2g_{1*}(\bar\Q_{Z_1})_3$ implies that $\tau$
is induced by an endomorphism of $R^1g_{0*}(\bar\Q_{Z_0})_3|_{\{u\}}$ of type $(-1,1)$.
Those form a $2$ dimensional vector space, corresponding to the deformations along $M_{5,1}$.

So a deformation $M_{5,2}\times N$, which does not lie in
$M_{5,2}\times M_{5,1}$, induces a non-zero flat section $\tau$ of
${\rm End}(\V\oplus\W)|_{M_{5,2}} $ of type $(-1,1)$.
which does not respect the direct sum decomposition. So one finds
a non-trivial morphism
$$
\tau: \V|_{M_{5,2}}= \bigoplus  R^2g_{1*}(\bar\Q_{Z_1})_3 \oplus
\bigoplus  R^2g_{1*}(\bar\Q_{Z_1})_2  \>>> \W|_{M_{5,2}}.$$
On the other hand,  $\W|_{M_{5,2}}$ is a direct sum
of several copies of the local systems
$$
R^2g_{1*}(\bar\Q_{Z_1})_1,\, R^2g_{1*}(\bar\Q_{Z_1})_4, \mbox{ \ \ and \ \ }
R^1f_*\bar \Q_X(1).
$$
Remark that the uniformization local system $R^2g_{1*}(\bar\Q_{Z_1})_2$ for $M_{5,2}$
is irreducible, as well as its dual $R^2g_{1*}(\bar\Q_{Z_1})_3$.

The local system $R^1f_*\bar\Q_X(1)$ is the variation of Hodge structures attached to the
universal family of plane curves of degree $5$, hence it is irreducible
by Deligne's irreducibility theorem \cite{DK}.
$R^2g_{1*}(\bar\Q_{Z_1})_1$ and $R^2g_{1*}(\bar\Q_{Z_1})_4$ are both irreducible, by a generalization
of Deligne's irreducibility theorem proved in \cite{VZ}, Lemma 4.1.

On the other hand, all the irreducible local systems considered above have different
Hodge types. So there exists no non-zero morphism between them, a contradiction.
\end{proof}
\begin{remark}
In \cite{VZ} we consider the subscheme
$$ M_{5,1}\times M_{5,1} \subset    M_{5,2}\times M_{5,1},$$
and the restriction of
$$h_2: \sZ_2\>>> M_{5,2}\times M_{5,1}$$
to this subscheme. It is shown there, that the set of CM-points $y\in M_{5,1}\times M_{5,1}$
is dense in $ M_{5,1}\times M_{5,1}$, i.e. the set of points $y$ for which the Hodge structure
$H^3(h_2^{-1}(y),\Q)$ has complex multiplication.

Since $M^s_{5,2}$ is a non-arithmetic ball quotient one should expect, according to the
Andr\'e-Oort conjecture, that the only positive dimensional component of
Zariski closure of the set of CM-points in $M_{5,2}\times M_{5,1}$
is $M_{5,1}\times M_{5,1}$.
\end{remark}

\end{document}